\documentclass[11pt]{article}
\usepackage{graphicx}
\usepackage{amsmath,amsfonts,amssymb}
\usepackage{url}
\usepackage{amsmath}
\usepackage{epsfig}
\usepackage{epstopdf}
\usepackage{color}
\usepackage{float}
\usepackage{subcaption}
\def\tto{\;{\lower 1pt \hbox{$\rightarrow$}}\kern -10pt
\hbox{\raise 2pt \hbox{$\rightarrow$}}\;}

\def\ra{\rangle}
\def\la{\langle}

\def\epsilon{\varepsilon}
\def\B{\Bbb B}

\def\R{\Bbb R}

\def\N{\Bbb N}
\def\ox{\bar{x}}

\def\ph{\varphi}

\newcounter{lk}

\begin{document}
\begin{center}
{\bf NESTEROV'S SMOOTHING TECHNIQUE AND MINIMIZING DIFFERENCES OF CONVEX FUNCTIONS FOR HIERARCHICAL CLUSTERING }\\[2ex]
 N. M. NAM\footnote{Fariborz Maseeh Department of Mathematics and Statistics, Portland State University, Portland, OR 97207, USA (mau.nam.nguyen@pdx.edu). Research of this author was partly supported by the National Science Foundation under grant \#1411817.}, W. GEREMEW\footnote{School of General Studies, Stockton University, Galloway, NJ 08205, USA (wondi.geremew@stockton.edu)}, S. REYNOLDS\footnote{Fariborz Maseeh Department of Mathematics and Statistics, Portland State University, Portland, OR 97207, USA (ser6@pdx.edu ).} and T. TRAN\footnote{Fariborz Maseeh Department of Mathematics and Statistics, Portland State University, Portland, OR 97207, USA (tuyen2@pdx.edu).}.\\[2ex]
\end{center}
\small{\bf Abstract.} A bilevel hierarchical clustering model is commonly used in designing  optimal multicast networks. In
this paper, we consider two different formulations of the bilevel hierarchical clustering
problem, a discrete optimization problem which can be shown to be NP-hard. Our approach is to reformulate
the problem as a continuous optimization problem by making some relaxations on the discreteness conditions.
Then Nesterov's smoothing technique and a numerical algorithm for minimizing differences of convex functions called the DCA are applied to cope with the nonsmoothness and nonconvexity of the problem. Numerical examples are provided to illustrate our method.

\newtheorem{Theorem}{Theorem}[section]
\newtheorem{Proposition}[Theorem]{Proposition}
\newtheorem{Remark}[Theorem]{Remark}
\newtheorem{Lemma}[Theorem]{Lemma}
\newtheorem{Corollary}[Theorem]{Corollary}
\newtheorem{Definition}[Theorem]{Definition}
\newtheorem{Example}[Theorem]{Example}
\renewcommand{\theequation}{\thesection.\arabic{equation}}
\normalsize
\vspace*{-0.2in}

\section{Introduction}
\vspace*{-0.1in}

Although convex optimization techniques and numerical algorithms have been the topics of extensive research for more than 50 years, solving large-scale optimization problems without the presence of convexity remains a challenge.  This is a motivation to search for new optimization methods that are capable of handling broader classes of functions and sets where convexity is not assumed. One of the most successful approaches to go beyond convexity is to consider the class of functions representable as  differences of two convex functions. Functions of this type are called \emph{DC functions}, where DC stands for \emph{difference of convex}. It was recognized early by P. Hartman \cite{H59} that the class of DC functions has many nice algebraic properties. For instance, this class of functions is closed under many operations usually considered in optimization such as taking linear combination, maximum, or product of a finite number of DC functions.

Given a linear space $X$, a \emph{DC program} is an optimization problem in which the objective function $f\colon X\to \R$  can be represented as $f=g-h$, where $g, h\colon X\to \R$ are convex functions.  This extension of convex programming enables us to take advantage of the available tools from convex analysis and optimization.  At the same time, DC programming is sufficiently broad to use in solving many nonconvex optimization problems faced in recent applications. Another feature of DC programming is that it possesses a very nice duality theory; see \cite{TA1} and the references therein.  Although DC programming had been known to be important for many applications much earlier, the first algorithm for minimizing differences of convex functions called  the \emph{DCA} was introduced by Tao and An in \cite{TA1,TA2}. The DCA  is a simple but effective optimization scheme used extensively in DC programming and its applications.

\emph{Cluster analysis} or \emph{clustering} is one of the most important problems  in many  fields such as  machine learning, pattern recognition, image analysis, data compression, and computer graphics.
Given a finite number of data points in a metric space, a centroid-based clustering problem seeks a finite number of cluster centers with each data point assigned to the nearest cluster center in a way that a certain distance, a measure of dissimilarity among data points, is minimized. Since many kinds of data encountered in practical applications have nested structures, they are required to use multilevel hierarchical clustering which involves grouping a data set into a hierarchy of clusters. In this paper, we apply the mathematical optimization approach to the bilevel hierarchical clustering problem.
In fact, using mathematical optimization in clustering is a very promising approach to overcome many disadvantages of the \emph{$k-$mean algorithm} commonly used in clustering; see \cite{abt,amt2,LBU03,nbg} and the references therein. In particular, the DCA was successfully applied in \cite{amt2} to a bilevel hierarchical clustering problem in which the distance measurement is defined by the squared Euclidean distance. Although the DCA in \cite{amt2} provides an effective way to solve the bilevel hierarchical clustering in high dimensions, it has not been used to solve the original model defined by the Euclidean distance measurement proposed in \cite{LBU03} which is  \emph{not suitable for the resulting DCA} according to the authors of \cite{amt2}.  By applying Nesterov's smoothing technique and the DCA, we are able to solve the original model proposed in \cite{LBU03} in high dimensions.


The paper is organized as follows. In Section 2, we present basic definitions and tools of optimization that are used throughout the paper. In section 3, we study two models of bilevel hierarchical clustering problems, along with two new algorithms based on Nesterov's smoothing technique and the DCA. Numerical examples and conclusions are presented in Section 4 and Section 5, respectively.

\section{Basic Definitions and Tools of Optimization}
\setcounter{equation}{0}
In this section, we  present two main tools of optimization used to solve the bilevel hierarchical crusting problem: the DCA introduced by Pham Dinh Tao and Nesterov's smoothing technique.

We consider throughout the paper DC programming:
\begin{equation}\label{dcf}
\mbox{\rm minimize}\, f(x):=g(x)-h(x), x\in \R^n,
\end{equation}
where $g\colon\R^n\to \R$ and $h\colon \R^n\to \R$ are convex functions. The function $f$ in \eqref{dcf} is called a \emph{DC function} and $g-h$ is called a \emph{DC decomposition} of $f$.

Given a convex function $g\colon \R^n\to \R$, the \emph{Fenchel conjugate} of $g$ is defined by
\begin{equation*}
g^*(y):=\sup\{\la y, x\ra -g(x)\; |\; x\in \R^n\}, y\in \R^n.
\end{equation*}
Note that  $g^*\colon \R^n\to (-\infty, \infty]$ is also a convex function. In addition,  $x\in \partial g^*(y)$ if and only if $y\in \partial g(x)$, where $\partial$ denotes the subdifferential operator in the sense of convex analysis; see, e.g., [19-21]

Let us present below the DCA introduced by Tao and An \cite{TA1,TA2} as applied to (\ref{dcf}). Although the algorithm is used for nonconvex optimization problems, the convexity of the functions involved still plays a crucial role.

{\bf The DCA}
\begin{center}
\begin{tabular}{| l |}
\hline
{\small INPUT}: $x_0\in \R^n$, $N\in \N$.\\
{\bf for} $k=1, \ldots, N$ {\bf do}\\
\qquad Find $y_k\in \partial h(x_{k-1})$.\\
\qquad Find $x_{k}\in \partial g^*(y_k)$.\\
{\bf end for}\\
{\small OUTPUT}: $x_{N}$.\\
\hline
\end{tabular}
\end{center}

Let us discuss below a convergence result of DC programming. A  function $h\colon \R^n\to \R$ is called $\gamma$-convex ($\gamma\geq 0$) if the function defined by $k(x):=h(x)-\frac{\gamma}{2}\|x\|^2$, $x\in \R^n$, is convex. If there exists $\gamma>0$ such that $h$ is $\gamma-$convex, then $h$ is called strongly convex.  We say that an element $\ox\in \R^n$ is a \emph{critical point} of the function $f$ defined by \eqref{dcf} if
\begin{equation*}\partial g(\ox)\cap \partial h(\ox)\neq \emptyset.
\end{equation*}
Obviously, in the case where both $g$ and $h$ are differentiable, $\ox$ is a critical point of $f$ if and only if $\ox$ satisfies the Fermat rule $\nabla f(\ox)=0$. The theorem below provides a convergence result for the DCA. It can be derived directly from \cite[Theorem 3.7]{TA2}.

\begin{Theorem} Consider the function $f$ defined by \eqref{dcf} and the sequence $\{x_k\}$ generated by the DCA. The following properties are valid:\\[1ex]
{\rm\bf(i)} If $g$ is $\gamma_1$-convex and $h$ is $\gamma_2$-convex, then
\begin{equation*}
f(x_{k})-f(x_{k+1})\geq \frac{\gamma_1+\gamma_2}{2}\|x_{k+1}-x_k\|^2\; \mbox{\rm for all }k\in \N.
\end{equation*}
{\rm\bf (ii)} The sequence $\{f(x_k)\}$ is monotone decreasing.\\[1ex]
 {\rm\bf (iii)} If $f$ is bounded from below,   $g$ is $\gamma_1$-convex and $h$ is $\gamma_2$-convex with $\gamma_1+\gamma_2>0$, and  $\{x_k\}$ is bounded, then every subsequential limit of the sequence $\{x_k\}$ is a critical point of $f$.
\end{Theorem}

Now we present a direct consequence of Nesterov's smoothing technique given in \cite{n}. In the proposition below, $d(x; \Omega)$ denotes the Euclidean distance and $P(x; \Omega)$ denotes the Euclidean projection from a point $x$ to a nonempty closed convex set $\Omega$ in $\R^n$.
\begin{Proposition}\label{p1} Given any $a\in \R^n$ and $\mu>0$, a Nesterov smoothing approximation of $\ph(x):=\|x-a\|$ defined in $\R^n$ has the representation
\begin{equation}\label{NA}
\ph_\mu(x):=\frac{1}{2\mu}\|x-a\|^2-\frac{\mu}{2}\big[d(\frac{x-a}{\mu}; \B)\big]^2.
\end{equation}
Moreover, $\nabla \ph_\mu(x)=P(\frac{x-a}{\mu}; \B)$ and
\begin{equation*}
\ph_\mu(x)\leq \ph(x)\leq \ph_\mu(x)+\frac{\mu}{2},
\end{equation*}
where $\B$ is the closed unit ball of $\R^n$.
\end{Proposition}

\section{The Bilevel Hierarchical Clustering Problem}\label{S:HC}
\setcounter{equation}{0}

Given a set of $m$ points (nodes) $a^1, a^2, \ldots, a^m$ in $\R^n$, our goal is to decompose this set into $k$ clusters. In each cluster, we would like to find a point $x^i$ among the nodes and assign it as the center for this cluster with all points in the cluster connected to this center. Then we will  find a \emph{total center} $x^*$ among the given points $a^1, a^2, \ldots, a^m$, and all centers are connected to this total center. The goal is to minimize the \emph{total transportation cost} in this tree computed by the sum of the distances from the total center to each center and from each center to the nodes in each cluster. This is a discrete optimization problem which can be shown to be \emph{NP-hard}. We will solve this problem based on continuous optimization techniques.

\subsection{The Bilevel Hierarchical Clustering: Model I}\label{S:HC}

The difficulty in solving this hierarchical clustering problem lies in the fact that the centers and total center have to be among the nodes. We first relax this condition with the use of artificial centers $x^1, x^2, \ldots, x^k$ that could be anywhere in $\R^n$. Then we define the total center as the centroid of  $x^1, x^2, \ldots, x^k$  given by
\begin{equation*}
x^*:=\frac{1}{k}(x^1+x^2+\cdots+x^k).
\end{equation*}
The total cost of the tree is given by
\begin{equation*}
\ph(x^1, \ldots, x^k):=\sum_{i=1}^m \min_{\ell =1, \ldots, k} \|x^\ell-a^i\| +\sum_{\ell=1}^k\|x^\ell-x^*\|.
\end{equation*}
Note that here each $a^i$ is assigned to its closest center. However, what we expect are the real centers, which can be approximated by trying to minimize the difference between the artificial centers and the real centers. To achieve this goal, define the function
\begin{equation*}
\phi(x^1, \ldots, x^k):=\sum_{\ell=1}^k \min_{i=1, \ldots,m}\|x^\ell-a^i\|, x^1, \ldots, x^k\in \R^n.
\end{equation*}
Observe that $\phi(x^1, \ldots, x^k)=0$ if and only if for every $\ell=1, \ldots, k$, there exists $i\in \{1, \ldots, m\}$ such that $x^\ell=a^i$, which means that $x^\ell$ is a real node. Therefore, we consider the constrained minimization problem:
\begin{align*}
\mbox{\rm minimize } \ph(x^1, \ldots, x^k) \;\mbox{\rm subject to } \phi(x^1, \ldots, x^k)=0.
\end{align*}
This problem can be converted to an unconstrained minimization problem:
\begin{equation}\label{HCP1}
\mbox{\rm minimize } f_{\lambda}(x^1, \ldots, x^k):=\ph(x^1, \ldots, x^k)+\lambda \phi(x^1, \ldots, x^k), x^1, \ldots, x^k\in \R^n,
\end{equation}
where $\lambda>0$ is a penalty parameter. Similar to the situation with the clustering problem, this new problem is nonsmooth and nonconvex, which can be solved by smoothing techniques and the DCA. Note that a particular case of this model in two dimensions was considered in \cite{LBU03} where the problem was solved using the derivative-free discrete gradient method established in \cite{Ba99}, but this method is not suitable for large-scale settings in high dimensions. The DCA was used in \cite{amt2} to solve a similar model with the squared Euclidean distance function used as the distance measurement. In this paper, the authors also addressed the difficulty of dealing with \eqref{HCP1} using the DCA as \emph{not suitable for the resulting
DCA.} Nevertheless, we will show in what follows that  the DCA is applicable to this model when combined with Nesterov's smoothing technique.

Note that the functions $\ph$ and $\phi$ in \eqref{HCP1} belong to the class of DC functions with the following DC decompositions:
\begin{align*}
\ph(x^1, \ldots, x^k)&=\sum_{i=1}^m\Big[ \sum_{\ell=1}^k\|x^\ell-a^i\|-\max_{r=1, \ldots,k}\sum_{\ell=1, \ell\neq r}^k\|x^\ell-a^i\|\Big]+\sum_{\ell=1}^k\|x^\ell-x^*\|\\
&=\Big[\sum_{i=1}^m\sum_{\ell=1}^k\|x^\ell-a^i\|+\sum_{\ell=1}^k\|x^\ell-x^*\|\Big]-\sum_{i=1}^m\max_{r=1, \ldots,k}\sum_{\ell=1, \ell\neq r}^k\|x^\ell-a^i\|,\\
\phi(x^1, \ldots, x^k)&=\sum_{\ell=1}^k\sum_{i=1}^m\|x^\ell-a^i\|-\sum_{\ell=1}^k\max_{s=1,\ldots,m}\sum_{i=1, i\neq s}^m\|x^\ell-a^i\|.
\end{align*}
It follows that the objective function $f_{\lambda}$ in  \eqref{HCP1} has the DC decomposition:
\begin{align*}
f_{\lambda}(x^1, \ldots, x^k)&=\Big[(1+\lambda)\sum_{\ell=1}^k\sum_{i=1}^m\|x^\ell-a^i\|+\sum_{\ell=1}^k\|x^\ell-x^*\|\Big]\\
&-\Big[\sum_{i=1}^m \max_{r=1, \ldots,k}\sum_{\ell=1, \ell\neq r}^k\|x^\ell-a^i\|+\lambda \sum_{\ell=1}^k\max_{s=1,\ldots,m}\sum_{i=1, i\neq s}^m\|x^\ell-a^i\|\Big].
\end{align*}
This DC decomposition is not suitable for applying the DCA because there is no closed form for a subgradient of the function $g^*$ involved.

In the next step, we apply Nesterov's smoothing technique from Proposition \ref{p1} to approximate the objective function $f_\lambda$ by a new DC function favorable for applying the DCA. To accomplish this goal, we simply replace each term of the form $\|x-a\|$ from the first part of $f_\lambda(x^1, \ldots, x^k)$ (the \emph{positive part}) by the smooth approximation \eqref{NA}, while keeping the second part (the \emph{negative part}) the same. As a result, we obtain
\begin{align*}
f_{\lambda\mu}(x^1, \ldots, x^k)&:=\; \frac{(1+\lambda)\mu}{2} \sum_{i=1}^{m}\sum_{\ell =1}^{k}\Bigg\|\frac{x^\ell - a^i}{\mu}\Bigg\|^2 \;   +   \; \frac{\mu}{2}\sum_{\ell =1}^{k}\Bigg\|\frac{x^\ell - x^*}{\mu}\Bigg\|^2 \;  \\
& -   \frac{(1+\lambda)\mu}{2} \sum_{i=1}^{m}\sum_{\ell =1}^{k}\left[d\left(\frac{x^\ell - a^i}{\mu};\mathbb B\right)\right]^2-\; \frac{\mu}{2}\sum_{\ell =1}^{k}\left[d\left(\frac{x^\ell - x^*}{\mu};\mathbb B\right)\right]^2 \\
&   -    \; \sum_{i=1}^{m}\max_{r=1, \ldots,k}\sum_{\ell=1, \ell\neq r}^k\|x^\ell - a^i\|\;   -    \;\lambda\sum_{\ell =1}^{k}\max_{s=1, \ldots,m}\sum_{i=1, i\neq s}^m\|x^\ell - a^i\|.
\end{align*}
The original bilevel hierarchical clustering problem now can be solved using a DC program:
\begin{align*}
\mbox{\rm minimize}\; f_{\lambda\mu}(x^1, \ldots, x^k)=g_{\lambda\mu}(x^1, \ldots, x^k) - h_{\lambda\mu}(x^1, \ldots, x^k),\;\;x^1, \ldots, x^k \in \mathbb{R}^ n.
\end{align*}
In this formulation, $g_{\lambda\mu}$ and $h_{\lambda\mu}$ are convex functions on $(\mathbb{R}^ n)^k$ defined by
\begin{align*}
&g_{\lambda\mu}(x^1, \ldots, x^k):= g^1_{\lambda\mu}(x^1, \ldots, x^k) + g^2_{\lambda\mu}(x^1, \ldots, x^k),\\
&h_{\lambda\mu}(x^1, \ldots, x^k):= h^1_{\lambda\mu}(x^1, \ldots, x^k) + h^2_{\lambda\mu}(x^1, \ldots, x^k) + h^3_{\lambda\mu}(x^1, \ldots, x^k) + h^4_{\lambda\mu}(x^1, \ldots, x^k),
\end{align*}
with their respective components defined as
\begin{align*}
&g^1_{\lambda\mu}(x^1, \ldots, x^k):=\frac{1+\lambda}{2\mu} \sum_{i=1}^{m}\sum_{\ell =1}^{k}\|x^\ell - a^i\|^2, \;  g^2_{\lambda\mu}(x^1, \ldots, x^k):= \frac{1}{2\mu}\sum_{\ell =1}^{k}\|x^\ell - x^*\|^2,\\
& h^1_{\lambda\mu}(x^1, \ldots, x^k):= \frac{(1+\lambda)\mu}{2} \sum_{i=1}^{m}\sum_{\ell =1}^{k}\left[d\left(\frac{x^\ell - a^i}{\mu};\mathbb B\right)\right]^2, \; h^2_{\lambda\mu}(x^1, \ldots, x^k):= \frac{\mu}{2}\sum_{\ell =1}^{k}\left[d\left(\frac{x^\ell - x^*}{\mu};\mathbb B\right)\right]^2,\\
& h^3_{\lambda\mu}(x^1, \ldots, x^k):= \sum_{i=1}^{m}\max_{r=1, \ldots,k}\sum_{\ell=1, \ell\neq r}^k\|x^\ell - a^i\|, \; 	h^4_{\lambda\mu}(x^1, \ldots, x^k):=\lambda\sum_{\ell =1}^{k}\max_{s=1, \ldots,m}\sum_{i=1, i\neq s}^m\|x^\ell - a^i\|.
\end{align*}

To facilitate the gradient and subgradient calculations for the DCA, we introduce a \emph{data matrix} $\mathbf A$ and a \emph{variable matrix} $\mathbf X$.  The data $\mathbf A$ is formed by putting each $a^i$, $i =1, \ldots, m$, in the $i^{th}$ row, i.e.,
\begin{equation*}
\mathbf{A} =
   \left(
   \begin{matrix}
      a_{11}      	&    a_{12}   	&    a_{13}		& \dots   		& a_{1n}\\
      a_{21}      	&    a_{22}   	&    a_{23}		& \dots   		& a_{2n}\\
      \vdots 		&  \vdots 		& \vdots  		& 		 	& \vdots\\
      a_{m1}      	&    a_{m2}   	&    a_{m3}	& \dots   		& a_{mn}
   \end{matrix}
   \right).
\end{equation*}
Similarly, if $x^1, \ldots, x^k$ are the $k$ cluster centers, then the variable $\mathbf X$ is formed by putting each $x^\ell$, $\ell =1, \ldots ,k$, in the $\ell^{th}$ row, i.e.,
\begin{equation*}
\mathbf{X} =
   \left(
   \begin{matrix}
      x_{11}      &    x_{12}   &    x_{13}& \dots   & x_{1n}\\
      x_{21}      &    x_{22}   &    x_{23}& \dots   & x_{2n}\\
      \vdots 		&  \vdots 		& \vdots  	& 	& \vdots\\
      x_{k1}      &    x_{k2}   &    x_{k3}& \dots   & x_{kn}
   \end{matrix}
   \right).
\end{equation*}
Then the variable matrix $\mathbf X$ of the optimization problem belongs to $\mathbb{R}^{k\times n}$, the linear space of $k$ by $n$ real matrices equipped with the inner product $\langle \mathbf X, \mathbf Y\rangle:=\mbox{\rm trace}(\mathbf X^T\mathbf Y)$. The \emph{Frobenius norm} on $\mathbb{R}^{k\times n}$ is defined by
\begin{equation*}
\big\|\mathbf X\big\|_F:= \sqrt{\Big\langle \mathbf X, \mathbf X \Big\rangle} =  \sqrt{\sum_{\ell =1}^{k}\langle x^\ell, x^\ell\rangle} = \sqrt{\sum_{\ell =1}^{k}\|x^\ell\|^2}.
\end{equation*}
Finally, we represent the average of the $k$ cluster centers by $x^*$, i.e., $x^*:=\frac{1}{k}\sum_{j=1}^{k}x^j$.

\subsection*{Gradient and Subgradient Calculations for the DCA}

Let us start by computing the gradient of
\begin{equation*}
g_{\lambda\mu}(\mathbf X) = g^1_{\lambda\mu}(\mathbf X) + g^2_{\lambda\mu}(\mathbf X).
\end{equation*}

Using the Frobenius norm, the function $g^1_{\lambda\mu}$  can equivalently be written as
\begin{align*}
g^1_{\lambda\mu}(\mathbf X) \;	&=\; \frac{1+\lambda}{2\mu} \sum_{i=1}^{m}\sum_{\ell =1}^{k}\|x^\ell - a^i\|^2\\
					&=\; \frac{1+\lambda}{2\mu} \sum_{i=1}^{m}\sum_{\ell =1}^{k}\left[\|x^\ell\|^2-2\langle x^\ell,  a^i\rangle  + \|a^i\|^2\right]\\
					&=\; \frac{1+\lambda}{2\mu} \left[m\big\|\mathbf X\big\|_F^2  -  2\Big\langle \mathbf X, \mathbf E_{km}\mathbf A\Big\rangle  +  k\big\|\mathbf A\big\|_F^2\right],
\end{align*}
where $\mathbf E_{km}$ is a $k\times m$ matrix whose entries are all ones. Hence, one can see that $g^1_{\lambda\mu}$ is differentiable and its gradient is given by
\[\nabla g^1_{\lambda\mu}(\mathbf X) \;	=\; \frac{1+\lambda}{\mu}\left[m\mathbf X - \mathbf E_{km}\mathbf A\right].\]															
Similarly, $g^2_{\lambda\mu}$ can equivalently be written as
\begin{align*}
g^2_{\lambda\mu}(\mathbf X) \; &=\; \frac{1}{2\mu}\sum_{\ell =1}^{k}\| x^\ell - x^*\|^2\\
					&=\; \frac{1}{2\mu} \sum_{\ell =1}^{k}\left[\|x^\ell\|^2-2\langle  x^\ell,  x^*\rangle  + \|x^*\|^2\right]\\
					&=\; \frac{1}{2\mu} \left[\big\|\mathbf X\big\|_F^2  -  \frac{2}{k}\Big\langle \mathbf X, \mathbf E_{kk}\mathbf X\Big\rangle  + \frac{1}{k}\Big\langle \mathbf X, \mathbf E_{kk}\mathbf X\Big\rangle\right]\\
					&=\; \frac{1}{2\mu} \left[\big\|\mathbf X\big\|_F^2  -  \frac{1}{k}\Big\langle \mathbf X, \mathbf E_{kk}\mathbf X\Big\rangle \right],
\end{align*}
where $\mathbf E_{kk}$ is a $k\times k$ matrix whose entries are all ones. Hence,  $g^2_{\lambda\mu}$ is differentiable and its gradient is given by
\[\nabla g^2_{\lambda\mu}(\mathbf X) \;	=\; \frac{1}{\mu}\left[\mathbf X - \frac{1}{k}\mathbf E_{kk}\mathbf X\right].\]
Since $g_{\lambda\mu}(\mathbf X) = g^1_{\lambda\mu}(\mathbf X) + g^2_{\lambda\mu}(\mathbf X)$, its gradient can be computed by
\begin{align*}
 \nabla g_{\lambda\mu}(\mathbf X)	&=\nabla g^1_{\lambda\mu}(\mathbf X)+\nabla g^2_{\lambda\mu}(\mathbf X)\\
						&=\frac{1+\lambda}{\mu}\left[m\mathbf X - \mathbf E_{km}\mathbf A\right]+\frac{1}{\mu}\left[\mathbf X - \frac{1}{k}\mathbf E_{kk}\mathbf X\right]\\
						&=\frac{1}{\mu}\left[(1+\lambda)m\mathbf X - (1+\lambda)\mathbf E_{km}\mathbf A+\mathbf X - \frac{1}{k}\mathbf E_{kk}\mathbf X\right]\\
						&=\frac{1}{\mu}\left[\left[\left[(1+\lambda)m + 1\right]\mathbf I_{kk} - \frac{1}{k}\mathbf E_{kk}\right]\mathbf X - (1+\lambda)\mathbf E_{km}\mathbf A\right].
\end{align*}
Therefore,
\[\nabla g_{\lambda\mu}(\mathbf X) = \frac{1}{\mu}\Big[\big(((1+\lambda)m + 1)\mathbf I_{kk} - \mathbf J\big)\mathbf X - (1+\lambda)\mathbf S\Big],\;\; \text{where}\;\; \mathbf J=\frac{1}{k}\mathbf E_{kk},\;\; \text{and}\;\; \mathbf S = \mathbf E_{km}\mathbf A.\]

Our goal now is to compute $\nabla g^*(Y)$, which can be accomplished by the relation
\begin{equation*}
\mathbf X=\nabla g^*(\mathbf Y)\; \mbox{\rm if and only if }\mathbf Y=\nabla g(\mathbf X).
\end{equation*}
The latter can equivalently be written as
\begin{align*}
&\left[\left[1+(1+\lambda)m\right]\mathbf I_{kk} - \mathbf J\right]\mathbf X 	= \left[ (1+\lambda)\mathbf S+\mu \mathbf Y\right].
\end{align*}
Then with some algebraic manipulation we can show that
\begin{equation*}
\nabla g^*(\mathbf Y)=\mathbf X = \left[\frac{1}{1+(1+\lambda)m}\mathbf I_{kk} \;+\; \frac{1}{[1+(1+\lambda)m](1+\lambda)m}\mathbf J\right]\left[ (1+\lambda)\mathbf S+\mu \mathbf Y\right].
\end{equation*}

Next, we will demonstrate in more detail the techniques we used in finding a subgradient for the convex function $h_{\lambda\mu}$. Recall that  $h_{\lambda\mu}$ is defined by
\begin{equation*}
h_{\lambda\mu}(\mathbf X) =\sum_{i=1}^4 h^i_{\lambda\mu}(\mathbf X).
\end{equation*}
We will start with the function $h^1_{\lambda\mu}$ given by
\[h^1_{\lambda\mu}(\mathbf X) \;=\; \frac{(1+\lambda)\mu}{2} \sum_{i=1}^{m}\sum_{\ell =1}^{k}\left[d\left(\frac{x^\ell - a^i}{\mu};\mathbb B\right)\right]^2.\]
From its representation, one can see that $h^1_{\lambda\mu}$ is differentiable, and hence its subgradient coincides with its gradient, that can be computed by  the partial derivatives with respect to $x^1, \cdots, x^k$, i.e.,
\begin{align*}
&\frac{\partial h^1_{\lambda\mu}}{\partial x^\ell }(\mathbf X) \; =\; (1+\lambda) \sum_{i=1}^{m}\left[\frac{x^\ell - a^i}{\mu} -P\left(\frac{x^\ell - a^i}{\mu}; \B \right)\right].
\end{align*}
Thus, for $\ell = 1, 2, \ldots, k$, $\nabla h^1_{\lambda\mu}(\mathbf X)$ is the $k\times n$ matrix $\mathbf U$ whose $\ell^{th}$ row is $\frac{\partial h^1_{\lambda\mu}}{\partial x^\ell }(\mathbf X)$.

Similarly, one can see that the function $h^2_{\lambda\mu}$ given by
\begin{equation*}h^2_{\lambda\mu}(\mathbf X) \; =\; \frac{\mu}{2}\sum_{\ell =1}^{k}\left[d\left(\frac{x^\ell -  x^*}{\mu};\mathbb {B}\right)\right]^2
\end{equation*}
is differentiable with its partial derivatives computed by
\begin{align*}
&\frac{\partial h^2_{\mu}}{\partial x^\ell }(\mathbf X) \; =\; \left[\frac{\mathbf x^\ell - x^*}{\mu} - P\left(\frac{x^\ell - x^*}{\mu}; \B \right)\right]\; -\;\frac{1}{k} \sum_{j=1}^{k}\left[\frac{x^j - x^*}{\mu} - P\left(\frac{ x^j -  x^*}{\mu}; \B \right)\right].&
\end{align*}
Hence, for $\ell = 1, 2, \ldots, k$, $\nabla h^2_{\mu}(\mathbf X)$ is the $k\times n$ matrix $\mathbf{V}$ whose $\ell^{th}$ row is $\frac{\partial h^2_{\mu}}{\partial x^\ell }(\mathbf X)$.

Unlike $h^1_{\lambda\mu}$ and $h^2_{\lambda\mu}$, the convex functions $h^3_{\lambda\mu}$ and $h^4_{\lambda\mu} $ are not differentiable, but both can be written as a finite sum of the maximum of a finite number of convex functions. Let us compute a subgradient of $h^3_{\lambda\mu}$ as an example. We have
\begin{equation*}
h^3_{\lambda\mu}(\mathbf X) \;=\; \sum_{i=1}^{m}\max_{r=1, \ldots,k}\sum_{\ell=1, \ell\neq r}^k\|x^\ell - a^i\|=\sum_{i=1}^m\gamma_i(\mathbf X),
\end{equation*}
where, for $i=1, \ldots, m$, \[\gamma_i(\mathbf X):=\max\left\{\gamma_{ir}(\mathbf X)=\sum_{\ell=1, \ell\neq r}^k\|x^\ell - a^i\|,\;\;r=1, \ldots,k\right\}.\]
Then, for each $i=1, \ldots, m$, we find $\mathbf W_i\in \partial \gamma_i(\mathbf X)$ according to the subdifferential rule for the maximum of convex functions. Then define $\mathbf W:=\sum_{i=1}^m\mathbf W_i$ to get a subgradient of the function $h^3_{\lambda\mu}$ at $\mathbf X$ by the subdifferential sum rule. To accomplish this goal, we first choose an index $r^*$ from the index set $\{1, \ldots, k\}$ such that \[\gamma_i(\mathbf X)=\gamma_{ir^*}(\mathbf X)=\sum_{\ell=1, \ell\neq r^*}^{k}\|x^\ell - a^i\|.\] Using the familiar subdifferential formula of the Euclidean norm function, the $\ell^{th}$ row $w_i^\ell$ for $\ell\neq r^*$ of the matrix $\mathbf{W}_i$ is determined as follows
\begin{equation*}
w_i^\ell:=\begin{cases}
 \frac{x^\ell - a^i}{\|x^\ell - a^i\|_2}        &\text{if}\; x^\ell \neq a^i,\\
\; u\in\mathbb{B}        &\text{if}\; x^\ell = a^i.
\end{cases}
\end{equation*}
The $r^{*th}$ row of the matrix $\mathbf W_i$ is $w_i^{r^*}:=0$.

The procedure for calculating a subgradient of the function $h^4_{\lambda\mu}$ given by
\begin{equation*}
h^4_{\lambda\mu}(x^1, \ldots, x^k) =\lambda\sum_{\ell =1}^{k}\max_{s=1, \ldots,m}\sum_{i=1, i\neq s}^m\|x^\ell - a^i\|,
\end{equation*}
 is very similar to what we just demonstrated for $h^3_{\lambda\mu}$.

At this point, we are ready to give a new DCA based algorithm for Model I.

{\bf Algorithm 1}.
\begin{center}
\begin{tabular}{| l |}
\hline
{\small INPUT}: $\mathbf X_0, \lambda_0,\mu_0, N \in \mathbb{N}$.\\
{\bf while} \mbox{stopping criteria} ($\lambda$, $\mu$) = false {\bf do}\\
\qquad {\bf for} $k = 1, 2, 3, \cdots, N\;\;$ {\bf do}\\
\qquad{\qquad Find $\mathbf Y_k\in \partial h_{1\lambda\mu}(\mathbf X_{k-1})$}.\\
\qquad{\qquad Find $\mathbf X_{k}\in \partial g_{1\lambda\mu}^*(\mathbf Y_k)$}.\\
\qquad {\bf end for}\\
\qquad update $\lambda \;\;\text{and} \; \;\mu$.\\
{\bf end while}\\
{\small OUTPUT}: $\mathbf X_{N}$.\\
\hline
\end{tabular}
\end{center}

\subsection{The Bilevel Hierarchical Clustering: Model II}\label{S:HC2}

In this section, we introduce the second model to solve the bilevel hierarchical clustering problem. In this model, we use an additional variable $x^{k+1}$ to denote the total center. At first we allow the total center $x^{k+1}$ to be a free point in $\mathbb{R}^n$, the same as the $k$ cluster centers. Then the total cost of the tree  is given by
\begin{equation*}
\ph(x^1, \ldots, x^{k+1}):=\sum_{i=1}^m \min_{\ell =1, \ldots, k} \|x^\ell-a^i\| +\sum_{\ell=1}^k\|x^\ell-x^{k+1}\|, x^1, \ldots, x^{k+1}\in \R^n.
\end{equation*}

To force the $k+1$ centers to be chosen from the given nodes (or to make them as close to the nodes as possible), we set the constraint
\begin{equation*}
\phi(x^1, \ldots, x^{k+1}):=\sum_{\ell=1}^{k+1} \min_{i=1, \ldots,m}\|x^\ell-a^i\| =0.
\end{equation*}
Our goal is to solve the optimization problem
\begin{align*}
&\mbox{\rm minimize } \ph(x^1, \ldots, x^{k+1})\\
&\mbox{\rm subject to } \phi(x^1, \ldots, x^{k+1}), x^1, \ldots, x^{k+1}\in \R^n.
\end{align*}
Similar to the first model, this problem formulation can be converted to an unconstrained minimization problem involving a penalty parameter $\lambda>0$:
\begin{equation}\label{HCP2}
\mbox{\rm minimize } f_{\lambda}(x^1, \ldots, x^{k+1}):=\ph(x^1, \ldots, x^{k+1})+\lambda \phi(x^1, \ldots, x^{k+1}), x^1, \ldots, x^{k+1}\in \R^n.
\end{equation}

Next, we apply Nesterov's smoothing technique to get an approximation of the objective function $f$ given in \eqref{HCP2} which involves two parameter $\lambda>0$ and $\mu>0$:
\begin{align*}
f_{\lambda\mu}(\mathbf X)&:=\; \frac{1+\lambda}{2\mu} \sum_{i=1}^{m}\sum_{\ell =1}^{k+1}\|x^\ell - a^i\|^2 \;  +  \;\frac{1}{2\mu}\sum_{\ell =1}^{k}\|x^\ell - x^{k+1}\|^2\; - \;\frac{1}{2\mu}\sum_{i=1}^{m}\|x^{k+1} - a^i\|^2 \\
 &-\; \frac{\lambda\mu}{2}\sum_{i=1}^{m}\left[d\left(\frac{x^{k+1} - a^i}{\mu};\mathbb {B}\right)\right]^2 -\; \frac{(1+\lambda)\mu}{2} \sum_{i=1}^{m}\sum_{\ell =1}^{k}\left[d\left(\frac{x^\ell - a^i}{\mu};\mathbb {B}\right)\right]^2\\
 &-\; \frac{\mu}{2}\sum_{\ell =1}^{k}\left[d\left(\frac{x^\ell - x^{k+1}}{\mu};\mathbb {B}\right)\right]^2\;-\;\sum_{i=1}^{m}\max_{r =1, \ldots, k}\sum_{\ell=1,  \ell\neq r}^{k}\|x^\ell - a^i\| \;-\;\lambda \sum_{\ell =1}^{k+1}\max_{s=1, \ldots, m}\sum_{i=1, i\neq s}^{m}\|x^\ell - a^i\|.
\end{align*}
As we will show in what follows, it is convenient to apply the DCA to  minimize the function $f_{\lambda\mu}$. This function can be represented as the differences of two convex functions defined on $\mathbb{R}^{(k+1)\times n}$ using a variable $X$ whose $i^{th}$ row is $x^i$ for $i=1, \ldots, k+1$:
\begin{equation*}
f_{\lambda\mu}(\mathbf X)=g_{\lambda\mu}(\mathbf X) - h_{\lambda\mu}(\mathbf X),\;\;\mathbf X \in \mathbb{R}^{(k+1)\times n}.
\end{equation*}
In this formulation, $g_{\lambda\mu}$ and $h_{\lambda\mu}$ are convex functions defined on $\mathbb{R}^{(k+1)\times n}$ by
\begin{equation*}
g_{\lambda\mu}(\mathbf X):= g^1_{\lambda\mu}(\mathbf X) + g^2_{\lambda\mu}(\mathbf X)
\end{equation*}
and
\begin{equation*}
h_{\lambda\mu}(\mathbf X):= h^1_{\lambda\mu}(\mathbf X) + h^2_{\lambda\mu}(\mathbf X) + h^3_{\lambda\mu}(\mathbf X) + h^4_{\lambda\mu}(\mathbf X)+ h^5_{\lambda\mu}(\mathbf X)+ h^6_{\lambda\mu}(\mathbf X),
\end{equation*}
with their respective components given by
\begin{align*}
& g^1_{\lambda\mu} (\mathbf X):= \frac{1+\lambda}{2\mu} \sum_{i=1}^{m}\sum_{\ell =1}^{k+1}\|x^\ell - a^i\|^2, \;   g^2_{\lambda\mu} (\mathbf X):= \frac{1}{2\mu}\sum_{\ell =1}^{k}\|x^\ell - x^{k+1}\|^2,\\
& h^1_{\lambda\mu} (\mathbf X):= \frac{1}{2\mu}\sum_{i=1}^{m}\|x^{k+1} - a^i\|^2, \;   h^2_{\lambda\mu} (\mathbf X):=\frac{\lambda\mu}{2}\sum_{i=1}^{m}\left[d\left(\frac{x^{k+1} - a^i}{\mu};\mathbb {B}\right)\right]^2, \\
& h^3_{\lambda\mu} (\mathbf X):=\frac{(1+\lambda)\mu}{2} \sum_{i=1}^{m}\sum_{\ell =1}^{k}\left[d\left(\frac{x^\ell - a^i}{\mu};\mathbb {B}\right)\right]^2, \;  h^4_{\lambda\mu} (\mathbf X):=\; \frac{\mu}{2}\sum_{\ell =1}^{k}\left[d\left(\frac{x^\ell - x^{k+1}}{\mu};\mathbb {B}\right)\right]^2,\\
&h^5_{\lambda\mu} (\mathbf X):=\sum_{i=1}^{m}\max_{r =1, \ldots, k}\sum_{\substack{\ell=1, \ell\neq r}}^{k}\|x^\ell - a^i\|, \;  h^6_{\lambda\mu} (\mathbf X):=\; \lambda\sum_{\ell =1}^{k+1}\max_{s=1, \ldots, m}\sum_{\substack{i=1, i\neq s}}^{m}\|x^\ell - a^i\|.
\end{align*}

\subsection*{Gradient and Subgradient Calculations for the DCA}

Let us start by computing the gradient of the first part  of the DC decomposition, i.e.,
\begin{equation*}
g_{\lambda\mu}(\mathbf X) = g^1_{\lambda\mu} (\mathbf X) + g^2_{\lambda\mu} (\mathbf X).
\end{equation*}
By applying similar techniques used in computing gradients/subgradients for Model I, $\nabla g^1_{\lambda\mu} (\mathbf X)$ can be written as
\begin{equation*}
\nabla g^1_{\lambda\mu} (\mathbf X) \;	=\; \frac{1+\lambda}{\mu}\left[m\mathbf X - \mathbf E\mathbf A\right],
\end{equation*}
where $\mathbf E$ is the $(k+1)\times m$ matrix whose entries are all ones. Similarly,  $\nabla g^2_{\lambda\mu} (\mathbf X)$ can also be written as
\[\nabla g^2_{\lambda\mu} (\mathbf X) \;	=\; \frac{1}{\mu}\left[\mathbf I + \mathbf T\right]\mathbf X,\]
where $\mathbf I$ is the $(k+1)\times (k + 1)$ identity matrix, and $\mathbf T$ is the  $(k+1)\times (k + 1)$ matrix whose entries are all zeros except its $(k+1)^{th}$ row and $(k+1)^{th}$ column, which both are filled by the vector $(-1, -1, -1, \ldots,-1, k-1)$.

It follows that
\begin{align*}
 \nabla g_{\lambda\mu}(\mathbf X)	&=\nabla g^1_{\lambda\mu} (\mathbf X) \;+\; \nabla g^2_{\lambda\mu} (\mathbf X)\\
						&=\frac{1+\lambda}{\mu}\left[m\mathbf X \;-\;\mathbf E\mathbf A \right] \;+\; \frac{1}{\mu}\left[\mathbf I \;+\;\mathbf T \right] \mathbf X\\
						&=\frac{1}{\mu}\left[c_1\mathbf I \;+\; \mathbf T \right]\mathbf X\;-\; \frac{1+\lambda}{\mu}\mathbf E\mathbf A,
\end{align*}
where $c_1 = 1 + (1+\lambda)m$. Our goal now is to compute $\nabla g_{\lambda\mu}^*(\mathbf Y)$, which can be accomplished by the relation
\begin{equation*}
\mathbf X=\nabla g_{\lambda\mu}^*(\mathbf Y)\; \mbox{\rm if and only if }\mathbf Y=\nabla g_{\lambda\mu}(\mathbf X).
\end{equation*}
The latter can be equivalently written as
\begin{align*}
\left[c_1\mathbf I \;+\; \mathbf T\right]\mathbf X 	&=  (1+\lambda)\mathbf E \mathbf A +\mu \mathbf Y,
\end{align*}
whose solutions can be explicitly computed by its $\ell^{th}$ row for $\ell=1, \ldots, k+1$:
\begin{align*}
&x^\ell = \frac{\left[  (1+\lambda)\mathbf E \mathbf A +\mu \mathbf Y\right]_\ell \; + \; x^{k+1}}{c_1}\; \mbox{\rm for }\ell = 1, \ldots, k,\\
&x^{k+1} = \frac{c_1\left[  (1+\lambda)\mathbf E \mathbf A +\mu \mathbf Y\right]_{k+1}\; + \; \sum_{\ell = 1}^k \left[  (1+\lambda)\mathbf E \mathbf A +\mu \mathbf Y\right]_\ell}{(c_1 + k) (c_1 - 1)}.
\end{align*}
In the representation
\begin{equation*}
h_{\lambda\mu}(\mathbf X) = h^1_{\lambda\mu}(\mathbf X) + h^2_{\lambda\mu}(\mathbf X) + h^3_{\lambda\mu}(\mathbf X) + h^4_{\lambda\mu}(\mathbf X)+ h^5_{\lambda\mu}(\mathbf X)+ h^6_{\lambda\mu}(\mathbf X),
\end{equation*}
the convex functions $h^1_{\lambda\mu} ,\,h^2_{\lambda\mu} ,\,h^3_{\lambda\mu},\,\text{and} \,h^4_{\lambda\mu}$ are differentiable. The  \emph{partial derivatives} of $h^1_{\lambda\mu}$ are given by
\begin{align*}
&\frac{\partial h^1_{\lambda\mu }}{\partial x^\ell}(\mathbf X) = 0\; \mbox{\rm for }\ell = 1, \ldots, k, \\
&\frac{\partial h^1_{\lambda\mu }}{\partial x^{k+1} }(\mathbf X)  = \frac{1}{\mu} \sum_{i=1}^{m}(x^{k+1} - a^i)  = \frac{1}{\mu}\left[ mx^{k+1} - \sum_{i=1}^m\mathbf A_i\right].
\end{align*}
Then $\nabla h^1_{\lambda\mu} (\mathbf X)$ is the $(k+1)\times n$ matrix $\mathbf{L}$ whose $\ell^{th}$ row is $\frac{\partial h^1_{\lambda\mu }}{\partial x^\ell }(\mathbf X)$  for $\ell = 1, \ldots,k+1$.
Similarly, the partial derivatives of $h^2_{\lambda\mu}$ are given by
\begin{align*}
&\frac{\partial h^2_{\lambda\mu }}{\partial x^\ell}(\mathbf X) = 0 \; \mbox{\rm for }\ell = 1, \ldots, k, \\
&\frac{\partial h^2_{\lambda\mu }}{\partial x^{k+1} }(\mathbf X)=\lambda \sum_{i=1}^{m}\left[\frac{x^{k+1} - a^i}{\mu}- P\left(\frac{x^{k+1} - a^i}{\mu}; \B \right)\right].
\end{align*}
Then $\nabla h^2_{\lambda\mu} (\mathbf X)$ is the $(k+1)\times n$ matrix $\mathbf{M}$ whose $\ell^{th}$ row is $\frac{\partial h^2_{\lambda\mu }}{\partial x^\ell }(\mathbf X)$,  for $\ell = 1, \ldots,k+1$.
Now we compute the gradient of $h^3_{\lambda\mu}$. We have

\begin{align*}
&\frac{\partial h^3_{\lambda\mu }}{\partial x^\ell }(\mathbf X) =(1+\lambda) \sum_{i=1}^{m}\left[\frac{x^\ell - a^i}{\mu} - P\left(\frac{x^\ell - a^i}{\mu}; \B \right)\right]\; \text{for}\;\ell = 1,\ldots, k,\\
&\frac{\partial h^3_{\lambda\mu }}{\partial x^{k+1} }(\mathbf X)=0.
\end{align*}
Thus, $\nabla h^3_{\lambda\mu} (\mathbf X)$ is the $(k+1)\times n$ matrix $\mathbf{U}$ whose $l^{th}$ row is $\frac{\partial h^3_{\lambda\mu }}{\partial x^\ell }(\mathbf X)$  for $\ell = 1, \ldots, ,k+1$.

Let us compute the gradient of $h^4_{\lambda\mu}$. We have
\begin{align*}
&\frac{\partial h^4_{\lambda\mu }}{\partial x^\ell }(\mathbf X) \; =\; \left[\frac{x^\ell - x^{k+1}}{\mu} - P\left(\frac{x^\ell - x^{k+1}}{\mu}; \B \right)\right]\;  \text{for}\; \ell = 1, \ldots, k,\\
&\frac{\partial h^4_{\lambda\mu }}{\partial x^{k+1} }(\mathbf X)=-\sum_{\ell =1}^{k}\left[\frac{x^\ell - x^{k+1}}{\mu} - P\left(\frac{x^\ell - x^{k+1}}{\mu}; \B \right)\right].
\end{align*}
Similarly, $\nabla h^4_{\lambda\mu} (\mathbf X)$ is the $(k+1)\times n$ matrix $V$ whose $\ell^{th}$ row is filled with $\frac{\partial h^4_{\lambda\mu }}{\partial x^\ell }(\mathbf X)$, for $\ell =1,\ldots,k+1$.

The procedure for computing subgradients of the last two nondifferentiable components, $h^5_{\lambda\mu} $ and $h^6_{\lambda\mu}$, is similar to the procedure we described for computing subgradients of $h^3_{\lambda\mu}$ in  Model I.

The following is a DCA based algorithm for  Model II.

{\bf Algorithm 2}.
\begin{center}
\begin{tabular}{| l |}
\hline
{\small INPUT}: $\mathbf X_0, \lambda_0,\mu_0, N \in \mathbb{N}$.\\
{\bf while} \mbox{stopping criteria} ($\lambda$, $\mu$) = false {\bf do}\\
\qquad {\bf for} $k = 1, 2, 3, \cdots, N\;\;$ {\bf do}\\
\qquad{\qquad Find $\mathbf Y_k\in \partial h_{\lambda\mu}(X_{k-1})$}.\\
\qquad{\qquad Find $\mathbf X_{k}\in \partial g_{\lambda\mu}^*(Y_k)$}.\\
\qquad {\bf end for}\\
\qquad update $\lambda \;\;\text{and} \; \;\mu$.\\
{\bf end while}\\
{\small OUTPUT}: $\mathbf X_{N}$.\\
\hline
\end{tabular}
\end{center}

\section{Numerical Experiments}

\begin{figure}[H]
\begin{subfigure}{.33\textwidth}
  \centering
  \includegraphics[width=0.95\textwidth]{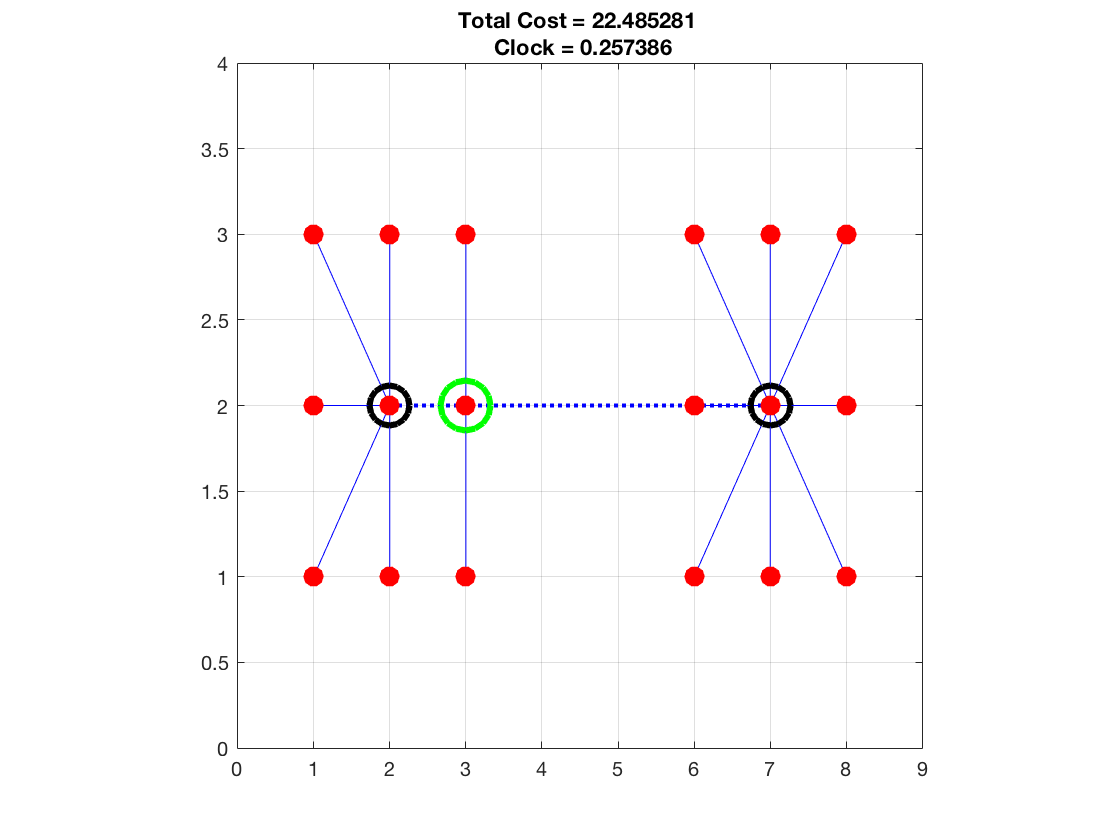}
  \caption{18 Data Points, 2 Centers}
  \label{fig:DS18}
\end{subfigure}%
\begin{subfigure}{.33\textwidth}
  \centering
\includegraphics[width=0.95\textwidth]{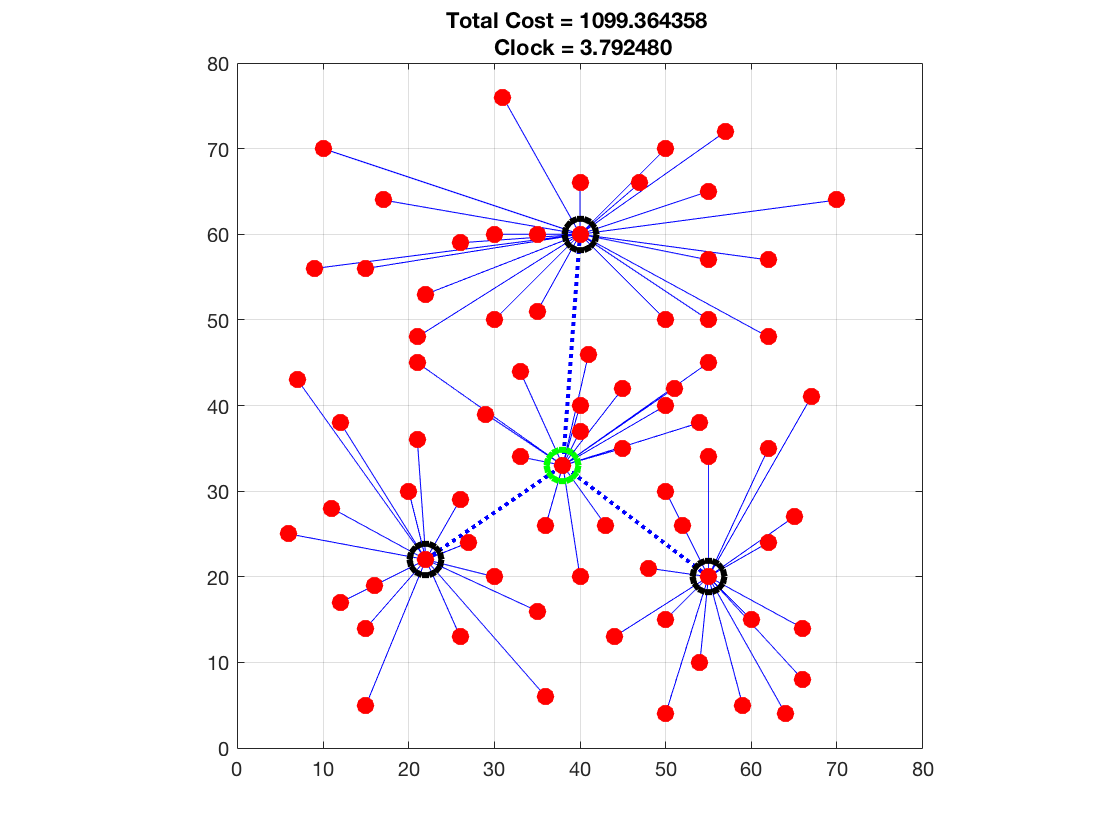}
  \caption{76 Data Points, 3 Centers}
  \label{fig:DS76}
\end{subfigure}
\begin{subfigure}{.33\textwidth}
  \centering
\includegraphics[width=0.95\textwidth]{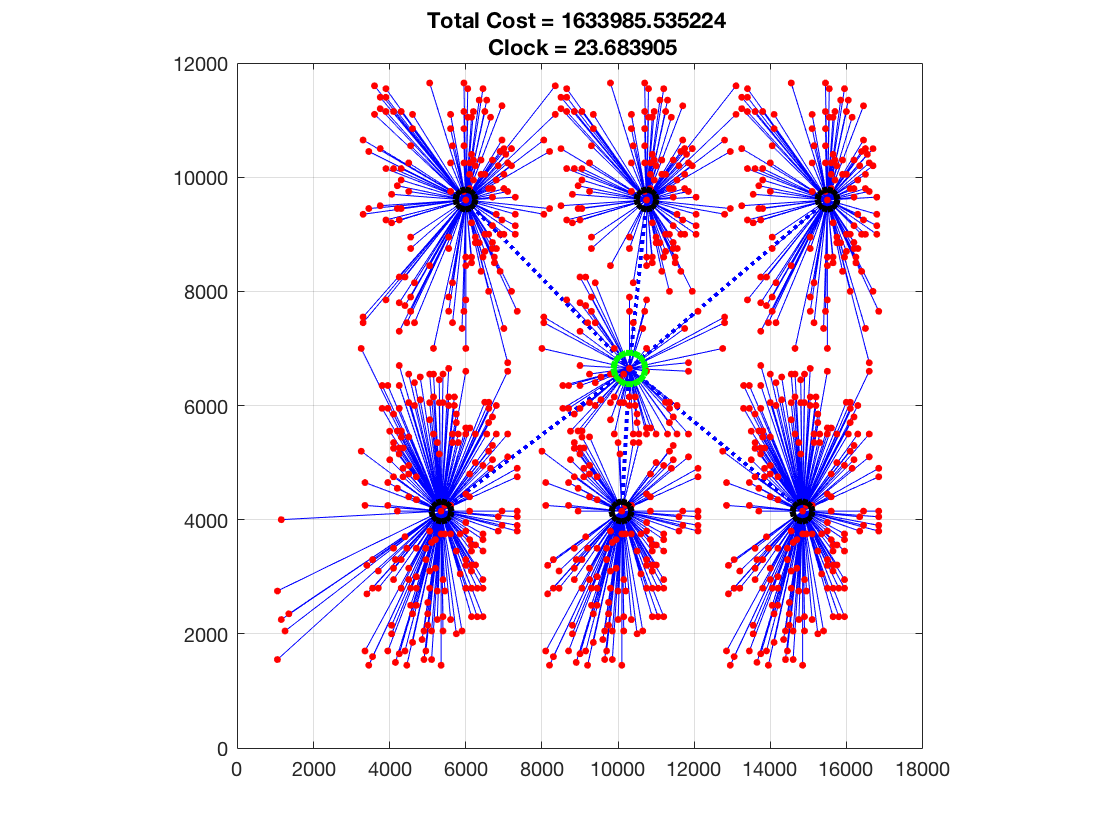}
  \caption{1002 Data Points, 6 Centers}
  \label{fig:DS1002}
\end{subfigure}
\caption{Plots of the three test data sets}
\label{fig:fig}
\end{figure}

We use MATLAB to code our algorithms and perform numerical experiments on a MacBook Pro with 2.2 GHz Intel Core i7 Processor, and 16 GB 1600 MHz DDR3 Memory. For our numerical experiments, we use three data sets: one artificial data set with 18 data points in $\R^2$ (see Figure \ref{fig:DS18}), the EIL76 and the PR1002 from  \cite{EIL-PR} (see Figure \ref{fig:DS76} and Figure \ref{fig:DS1002}), respectively.

The two MATLAB codes used to implement our two algorithms have two major parts: an \emph{outer} loop for updating the penalty and the smoothing parameters and an \emph{inner} loop for updating the cluster centers. The  penalty parameter $\lambda$ and the smoothing parameter $\mu$ are updated as follows. Choosing $\mu_0>0$, $\lambda_0>0$ and $\sigma_1>1$, $\sigma_2\in (0,1)$, we update $\lambda_{i+1}=\sigma_1\lambda_{i}$ and $\mu_{i+1}=\sigma_2\mu_{i}$ for $i\geq 0$ after each outer loop. To choose $\sigma_1$ and $\sigma_2$, we first let $N$ be the number of outer iterations and choose $\lambda_0,\lambda_\textup{max}$ as the initial and final values of $\lambda$, respectively, and similarly for $\mu_0,\mu_\textup{min}$. Then choose growth/decay parameters according to $\sigma_1=\left(\lambda_\textup{max}/\lambda_0\right)^{1/N}$ and $\sigma_2=\left(\mu_\textup{min}/\mu_0\right)^{1/N}$.

By trial and error we find that the values chosen for $\lambda_0$, $\lambda_\textup{max}$, $\mu_0$, and $\mu_\textup{min}$ in large part determine the performance of the two algorithms for each data set. Intuitively, we see that very large values of $\lambda$ will over-penalize the distance between an artificial center and its nearest data node and may prevent the algorithm from clustering properly. We therefore use $\lambda_0\leq1\ll\lambda_\textup{max}$ so that the algorithm has a chance to cluster the data before the penalty parameter takes effect. Similarly, we choose $\mu_\textup{min}\ll1\leq\mu_0$.

 We select the starting center $\mathbf X_0$ at a certain radius, $\gamma\; \text{rad}(\mathbf A)$, from the median point, $\text{median}(\mathbf A)$, of the entire data set, i.e.,
\[\mathbf X_0 = \text{median}(\mathbf A) + \gamma\; \text{rad}(\mathbf A)\;\mathbf{U},\] where $\text{rad}(\mathbf A) :=\max\{\|a^i - \text{median}(\mathbf A)\|\;|\; a^i\in\mathbf A\}$, $\gamma$ is a randomly chosen real number from the standard uniform distribution on the open interval (0,1), $\mathbf{U}$ is a $k \times n$ matrix whose $k$ rows are randomly generated unit vectors in $\mathbb{R}^n$, and the sum is in the sense of adding a vector to each row of a matrix.

\begin{figure}[H]
  \centering
\includegraphics[scale = 0.25]{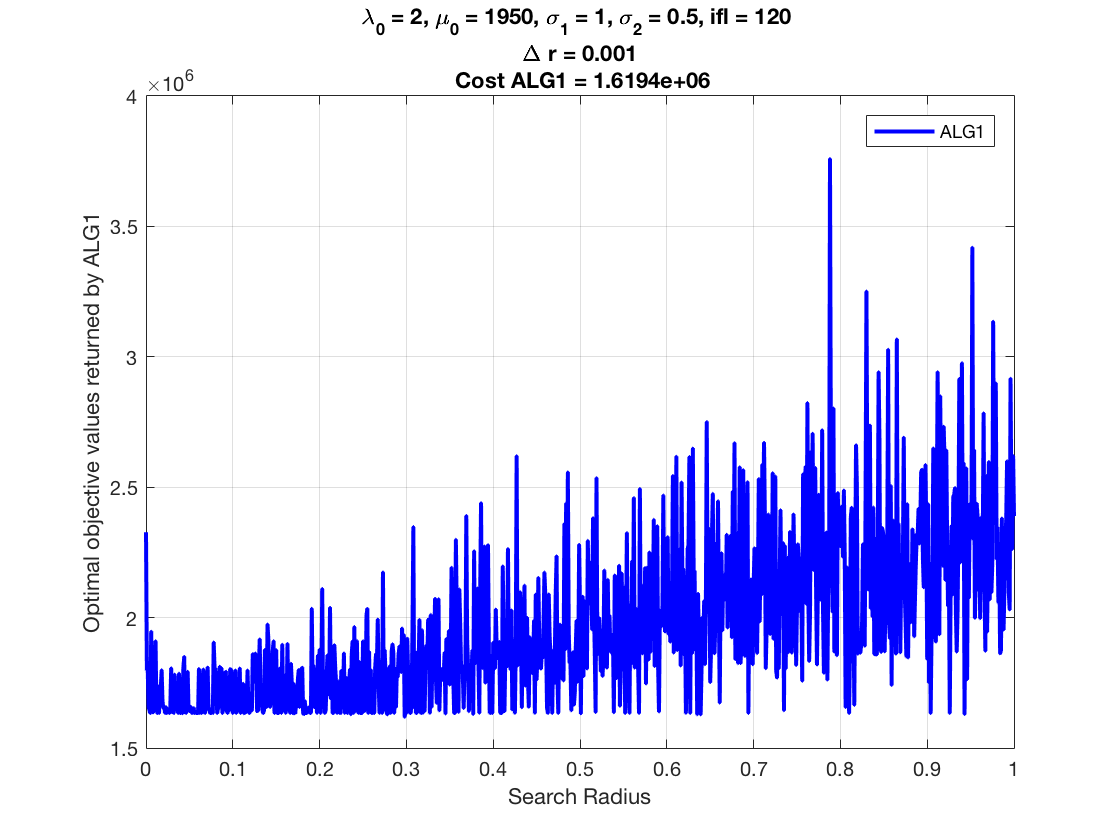}
\caption{PR1002, The 1002 City Problem, with 6 Cluster Centers}
\label{fig:rads}
\end{figure}

As showed in Tables \ref{tab:DS18R}, \ref{tab:EIL76R}, and \ref{tab:PR1002R}, both algorithms identify the optimal solutions with reasonable amount of time for both DS18 and EIL76 with two and three cluster centers, respectively. To get a good starting point which yields a better estimate of the optimal value for bigger data sets such as PR1002, we use a method called \emph{radial search} described as follows. Given initial radius $r_0>0$ and $m\in \N$, set $\gamma = i r_0$ for $i = 1, \ldots, m$. Then we test the algorithm with different starting points given by $\mathbf X_0(i) = \text{median}(\mathbf A) + i r_0(\text{rad}(\mathbf A)\mathbf{U})$ for $i = 1, \ldots, m$. Figure \ref{fig:rads} shows the result of the method applied to PR1002 with six cluster centers, where the $y$-axis represents the optimal value returned by ALG1 with different starting points $\mathbf X_0(i)$, as represented on the $x$-axis.
\begin{table}[H]
   \centering
   $\mu_0 = 5.70,\;\lambda_0 = 0.001,\; \sigma_1  = 7500,\;\sigma_2 = 0.5$ %
  \begin{tabular}{ c c c | c  c | c c | c  c c  }
\hline
&Cost1	&Cost2	&Time1	&Time2	&Iter1	&Iter2	&k	&m	&n	 \\ \hline
ADS18	&22.4853	&22.4853	&0.0690125	&0.0853712	&124	&124	&2	&18	&2	 \\
ADS18	&22.4853	&22.4853	&0.0678142	&0.0902829	&124	&124	&2	&18	&2	 \\
ADS18	&22.4853	&22.4853	&0.0694841	&0.0970158	&124	&124	&2	&18	&2	 \\
\hline
\end{tabular}
   \caption{Results for the 18 points artificial  data set.}
   \label{tab:DS18R}
\end{table}

\begin{table}[H]
   \centering
   $\mu_0 = 100,\;\lambda_0 = 10^{-6},\; \sigma_1  = 1,\;\sigma_2 = 0.5$
  \begin{tabular}{ c c c | c  c | c c | c  c c  }
\hline
&Cost1	&Cost2	&Time1	&Time2	&Iter1	&Iter2	&k	&m	&n	 \\ \hline
EIL76	&1125.48	&1107.47	&1.6126	&1.86261	&540	&540	&3	&76	&2	 \\
EIL76	&1099.36	&1099.36	&1.52001	&1.81243	&540	&540	&3	&76	&2	 \\
EIL76	&1099.36	&1099.36	&1.52914	&1.86536	&540	&540	&3	&76	&2	 \\
\hline
\end{tabular}
   \caption{Results for EIL76 data set.}
   \label{tab:EIL76R}
\end{table}

\begin{table}[H]
   \centering
    $\mu_0 = 1950,\;\lambda_0 = 10^{-6},\; \sigma_1  = 7500,\;\sigma_2 = 0.5$
  \begin{tabular}{ c c c | c  c | c c | c  c c  }
\hline
&Cost1	&Cost2	&Time1	&Time2	&Iter1	&Iter2	&k	&m	&n	 \\ \hline
PR1002	&1.63399e+06	&1.63399e+06	&22.1578	&25.3651	&330	&330	&6	&1002	&2	 \\
PR1002	&1.63399e+06	&1.63399e+06	&22.2978	&25.5373	&330	&330	&6	&1002	&2	 \\
PR1002	&1.63399e+06	&1.63399e+06	&23.5394	&26.0762	&330	&330	&6	&1002	&2	 \\
\hline
\end{tabular}
   \caption{Results for PR1002 data set.}
   \label{tab:PR1002R}
\end{table}

For comparison purposes, both Cost1 and Cost2 are computed by the same way. First, we systematically reassign the $k$ cluster centers returned by the respective algorithms by $k$ real nodes that are close to them, i.e., for $\ell = 1, \ldots, k$
\[ \bar{x}^\ell \;\;= {\operatorname{argmin}}\{\|x^\ell - a^i \|\;|\; \;a^i\in A\}. \]
Then the total center $x^*$ will be a real node, from the remaining nodes, whose sum of distances from the $k$ reassigned centers is the minimal, i.e.,
\[ x^*:= {\operatorname{argmin}}\Big\{\sum_{\ell =1}^{k}\|\bar{x}^\ell - a^i \|\;|\; a^i\in A \Big\}. \] The total cost is computed by adding the distance of each real node to its closest center (including the total center), and the distances of the total center from the $k$ cluster centers, i.e.,
\[\mbox{\rm Cost}:=\sum_{i=1}^m \min_{\ell =1, \ldots, k+1} \|\bar{x}^\ell-a^i\| +\sum_{\ell=1}^k\|\bar{x}^\ell-x^{k+1}\|,\;\text{where} \;x^{k+1} = x^*.\]

\section{Conclusions}
In our numerical experiments, one of the major challenge we face is how to choose optimal parameters for our algorithms. We observe that parameter selection is the decisive factor in terms of accuracy and speed of convergence of our proposed algorithms. The performance of the proposed algorithms  highly depends on the initial values set to the penalty and smoothing parameters, $\lambda_0$ and $\mu_0$; and their respective growth/decay factors $\sigma_1$ and $\sigma_2$. A better guidance of selecting the optimal parameter than the methods we suggest in this paper will be addressed in or future research.



\small
\begin{thebibliography}{99}

\bibitem{abt} L. T. H. An, M. T. Belghiti, P. D. Tao:  A new efficient algorithm based on DC programming and DCA for clustering. J. Glob. Optim.,   \textbf{27},  503--608 (2007)

 \bibitem{amt2}  L. T. H. An, and L. H. Minh, Optimization based DC programming and DCA
for hierarchical clustering. \textit{European J. Oper. Res}. \textbf{183} (2007), 1067--1085.


\bibitem{TA1} L. T. H. An, and P. D. Tao,  Convex analysis approach
to D.C. programming: Theory, algorithms and applications. \textit{Acta Math. Vietnam.} \textbf{22} (1997), 289--355.


\bibitem{Ba99} A. M. Bagirov, Derivative-free methods for
unconstrained nonsmooth optimization and its
numerical analysis. \textit{Investigacao Operacional}. \textbf{19} (1999), 75--93.


\bibitem{LBU03} A. Bagirov, Long Jia, I. Ouveysi, and A.M. Rubinov, Optimization based clustering algorithms in Multicast group hierarchies, in: Proceedings of the Australian Telecommunications, Networks and Applications Conference (ATNAC), 2003, Melbourne Australia
(published on CD, ISNB 0-646-42229-4).


\bibitem{bc} H. H. Bauschke and P. L. Combettes, {\em Convex Analysis and Monotone Operator Theory in Hilbert Spaces} Springer, New York, 2011.


\bibitem{bl} J. M. Borwein and A. S. Lewis, {\em Convex Analysis and Nonlinear Optimization}, 2nd edition, Springer, New York, 2006.


\bibitem{Ra} R. I. Bo\c t, {\em Conjugate Duality in Convex Optimization}, Springer, Berlin, 2010.


\bibitem{TA2} T. Pham Dinh and H. A. Le Thi, A d.c. optimization algorithm
for solving the trust-region subproblem, SIAM J. Optim. 8
(1998), 476--505.


\bibitem{H59} P. Hartman, On functions representable as a difference of convex functions.
\textit{Pacific J. Math}. \textbf{9}, (1959), 707--713.


\bibitem{HU} J.-B. Hiriart-Urruty, C. Lemar\'echal, {\em Convex Analysis and Minimization Algorithms I, II}, Springer, Berlin, 1993.


\bibitem{HU85} J. B. Hiriart-Urruty,  Generalized differentiability, duality and optimization for
problems dealing with differences of convex functions. \textit{Lecture Note in Economics and Math. Systems}. \textbf{256} (1985), 37--70.


\bibitem{m-book1} B. S. Mordukhovich, {\em Variational Analysis and Generalized Differentiation, I: Basic Theory, II: Applications}, Springer, Berlin, 2006.


\bibitem{bmn} B. S. Mordukhovich and N. M. Nam, {\em An Easy Path to Convex Analysis and Applications}, Morgan \& Claypool Publishers, San Rafael, CA, 2014.


\bibitem{nbg}N. M. Nam, R.B. Rector, D. Giles: Minimizing Differences of Convex Functions with Applications to Facility Location and Clustering, submitted.

\bibitem{n} Y. Nesterov:  Smooth minimization of non-smooth functions. Math. Program. \textbf{103}, 127--152 (2005)

\bibitem{EIL-PR} G. Reinelt, TSPLIB: A Traveling Salesman Problem Library. {\em ORSA Journal of Computing}. \textbf{3}
(1991), 376--384.


\bibitem{r} R. T. Rockafellar, {\em Convex Analysis}, Princeton University Press, Princeton, NJ, 1970.


\bibitem{r1} R. T. Rockafellar, {\em Conjugate Duality and Optimization}, SIAM, Philadelphia, PA, 1974.


\end{thebibliography}
\end{document}